\documentclass[12pt]{elsarticle}

\usepackage[top=99pt, bottom=75pt, left=72pt, right=70pt]{geometry}
\usepackage[utf8]{inputenc}
\usepackage{amsthm}
\usepackage{amsmath}
\usepackage{amssymb}
\usepackage{amscd,amsxtra}
\usepackage{MnSymbol}
\usepackage{mathtools}
\usepackage{tensor}
\usepackage{faktor}
\usepackage{dsfont}
\usepackage[all]{xy}
\usepackage{graphicx}
\usepackage[rgb]{xcolor}
\usepackage{tikz, tikz-cd}
\usetikzlibrary{arrows, matrix, intersections, math, calc, decorations.pathmorphing, decorations.markings, positioning}
\usepackage{enumitem}
\usepackage{hyperref}
\definecolor{R}{RGB}{255, 0, 34}
\definecolor{B}{RGB}{0, 85, 238}
\hypersetup{colorlinks = true,
            linkcolor = R,
            urlcolor  = B,
            citecolor = B,
            anchorcolor = B}

\newtheorem{theorem}{Theorem}

\newtheorem{proposal}{Proposal}

\numberwithin{theorem}{subsection}
\numberwithin{proposition}{subsection}
\numberwithin{lemma}{subsection}
\numberwithin{claim}{subsection}
\numberwithin{corollary}{subsection}
\numberwithin{conjecture}{subsection}
\numberwithin{definition}{subsection}
\numberwithin{remark}{subsection}
\newtheorem{example}[theorem]{Example}
\newtheorem{conjetura}{Conjecture}
\newtheorem{teorema}{Theorem}

\newcommand{\beq}{\begin{equation}}
\newcommand{\eeq}{\end{equation}}
\newcommand{\beqa}{\begin{eqnarray}}
\newcommand{\eeqa}{\end{eqnarray}}
\newcommand{\beaa}{\begin{eqnarray*}}
\newcommand{\ben}{\begin{eqnarray*}}
\newcommand{\eaa}{\end{eqnarray*}}
\newcommand{\een}{\end{eqnarray*}}

\makeatletter
\def\ps@pprintTitle{%
  \let\@oddhead\@empty
  \let\@evenhead\@empty
  \def\@oddfoot{\reset@font\hfil\thepage\hfil}
  \let\@evenfoot\@oddfoot
}

\journal{}

\begin{document}

\begin{frontmatter}



\title{\textbf{A tentative proposal towards an equivariant mirror symmetry for Hitchin systems}} 


\author{John Alexander Cruz Morales}

\ead{jacruzmo@unal.edu.co} 

\affiliation{organization={Universidad Nacional de Colombia, Departamento de Matemáticas},
            addressline={Ciudad Universitaria}, 
            city={Bogotá},
           country={Colombia}}

            
\begin{abstract}
\small{Motivated by Aganagic's equivariant mirror symmetry for certain Coulomb branches of a $3d$ $\mathcal{N}= 4$ gauge quiver theory, we would like to propose a set of ideas towards an extension of Aganagic's proposal to Hitchin systems. At the end, there are two main points in our proposal, namely that the \emph{equivariant mirror} of the Hitchin systems should be a Landau-Ginzburg model (with twisted masses) and that the duality between additive and multiplicative varieties in the context of mirror symmetry for Nakajima quiver varieties should be considered in the case of Hitchin systems.} 
\end{abstract}





\end{frontmatter}

\setlength{\parindent}{0cm}

\section{\textbf{Introduction}}

For more than three decades, mirror symmetry has captured the attention and interest of mathematicians and physicists, and a lot of work has been done. Different kinds of dualities in dimensions two, three and four are now packaged under the name \emph{mirror symmetry}. The usual one in mathematics is the two-dimensional mirror symmetry, as duality between symplectic and complex geometry. This appears in different flavours like topological, homological, Hodge theoretic or SYZ mirror symmetry. On the other hand, three-dimensional mirror symmetry is known in mathematics as symplectic duality and four-dimensional mirror symmetry is closely related to the geometric Langlands program. In this note, we will mainly be concerned with two-dimensional mirror symmetry; however, the main objects are certain holomorphic symplectic manifolds, which are also central in the three-dimensional and four-dimensional picture. Therefore, in this text, \emph{mirror symmetry} means two-dimensional mirror symmetry.  \\

It would be difficult to try to measure the impact that the work in mirror symmetry has had in both mathematics and physics. Many important questions related to the mirror phenomenon are still open and new problems have appeared in the landscape. The research in mirror symmetry nowadays is as strong as it was thirty years ago. The main aim of this text is to present a set of ideas that we believe could contribute to improving our understanding of what mirror symmetry is in a particular context. \\

The main motivation for our proposal is the interest in understanding mirror symmetry for noncompact holomorphic symplectic varieties with a torus action. In fact, we are interested in the case where those manifolds are not just holomorphic symplectic but hyperk\"ahler. The two main examples of such objects we have in mind are Nakajima quiver varieties and Hitchin systems\footnote{One of the original motivations of writing this note was to study the relations between Nakajima quiver varieties and Hitchin systems. Some work in this direction can be found in  \cite{schap}, where the so-called comet quiver varieties are studied. It would be interesting to see how \cite{schap} could be related to the discussion we are going to present here.}. Concretely, we are looking for a setting in which to discuss equivariant mirror symmetry for Hitchin systems. Therefore, some clarification is needed about what we mean by \emph{equivariant mirrors}. 

\subsection{\textbf{Equivariant mirror symmetry}}

In \cite{Teleman}, Teleman describes a framework to study equivariance in two-dimensional extended topological quantum field theories. As it is indicated in the introduction of loc.cit., this is a story of the categorified representation theory of a compact Lie group. Teleman's work fits in the setting of boundary conditions in three-dimensional topological quantum field theories associated with a holomorphic symplectic manifold $X$, known as Rozansky-Witten theory. His approach can be summarized in his own words as: \\

``\emph{Pure topological gauge theory in 3 dimensions for a compact Lie group $G$ is equivalent to the Rozansky-Witten theory for the BFM\footnote{This is a space introduced in \cite{bmf} which is closely related to the cotangent bundle to the space of conjugacy classes in $G^{\vee}_\mathbb{C}$.} space of the Langlands dual Lie group $G^{\vee}$}'' \\

Despite the relevance of Teleman's work, which might play an important role in the story we are trying to describe in this note, we approach the equivariance in a slightly different context\footnote{The literature for the equivariant approach to mirror symmetry is quite extensive. We do not pretend to be exhaustive. There are many relevant works in the equivariant setting closely related to the topics of this note. See, for instance, \cite{andersen, halpern, hausel} for discussion related to moduli of Higgs bundles, \cite{bom, MO} for discussion concerning equivariant quantum cohomology of Nakajima quiver varieties, and \cite{lekili} for discussion in the context of Fukaya categories. However, in a certain sense, these works only cover one half of Aganagic's proposal, as we will see.}. \\

In \cite{ag1, ag2, ag}, Aganagic proposed an equivariant mirror for the Coulomb branches $\mathcal{X}$ of certain $3d$ $\mathcal{N}= 4$ theories, which is a holomorphic symplectic variety. In a few words, the equivariant mirror of $\mathcal{X}$, in Aganagic's sense, is the ordinary mirror,  i.e., non-equivariant mirror, of its core. This could be a bit confusing, since the equivariant mirror turns out to be non-equivariant, but let us explain a little bit of Aganagic's setting. We will give a more detailed discussion in section \ref{agan}. \\

Starting with $\mathcal{X}$ of dimension $2m$, Aganagic constructs its two-dimensional mirror $\mathcal{Y}$ in the usual sense. We can think of $\mathcal{Y}$ as a \emph{multiplicative} version of $\mathcal{X}$. Since there is a torus acting on $\mathcal{X}$, we can think of $\mathcal{Y}$ as the equivariant mirror in the more standard way. However, Aganagic does not call it the equivariant mirror; instead, she calls it the \emph{upstairs mirror}. On the other hand, roughly speaking, the \emph{core} $X$ is the locus in $\mathcal{X}$  preserved by the symmetry that scales its holomorphic symplectic form; in particular, $X$ sits inside $\mathcal{X}$ as a holomorphic Lagrangian, in particular, the dimension of $X$ is $m$.\\

Now, one can construct the ordinary mirror $Y$ of $X$. This is the \emph{downstairs mirror}. It turns out that $\mathcal{Y}$ fibers over $Y$ with holomorphic Lagrangian $(\mathbb{C}^{\times})^m$ fibers. The main point in this framework is that working equivariantly, $X$, and its mirror $Y$, have as much information about the geometry as does $\mathcal{X}$. In this sense, Aganagic calls $Y$ the \emph{equivariant mirror} of $\mathcal{X}$\footnote{We might call it the \emph{up-down equivariant mirror} to avoid confusion with the more standard use of the word equivariant in the literature. However, at least for this note, we will keep Aganaigic's terminology.}. In this construction, some potentials $\mathcal{W}$ and $W$ appear, so in fact the mirror partners of $\mathcal{X}$ and $X$ are Landau-Ginzburg models. This is an important fact that will be discussed in section \ref{agan}. \\

The previous discussion can be summarized in the following diagram

\[
\begin{tikzcd}[column sep=6em, row sep=4em]
\mathcal{X} 
  \arrow[r, rightarrow, "upstairs-mirror"] 
  \arrow[dr, "equiv.mirror" description] 
& \mathcal{Y} \arrow[d, rightarrow] \\
X \arrow[r, rightarrow, "downstairs-mirror"']
 \arrow[u, hook]  
& Y
\end{tikzcd}
\]

The main idea of this text is to argue how the equivariant mirror symmetry described above, for Coulomb branches of a $3d$ $\mathcal{N}=4$ theories, could be implemented for the case of Hitchin systems. To do that, following some ideas from physics, we want to focus our attention on the Hitchin systems associated with the so-called class $\mathcal{S}$ theories. See section \ref{classS} for the discussion of such theories. Due to many formal similarities between these two contexts, the idea of extending Aganagic's picture to Hitchin systems seems very natural and plausible. \\ \\

The paper is organized as follows. The first three sections can be considered a short survey of some of the mathematical and physical topics which are relevant for the study of equivariant mirror symmetry for Hitchin systems. This is not an exhaustive survey and just represents our understanding of the subject. The last section is wildly speculative and presents our \emph{tentative proposal}, i.e., a set of ideas based on speculations and expectations on how the equivariant mirror symmetry, in Aganagic's sense, can work out for Hitchin systems. We do not formulate our proposals as conjectures since we do not have enough evidence for them. Only some formal analogies, similarities, and a certain feeling of structural beauty have guided our ideas.\\

Overall, there are two main messages in this note that we would like to emphasize, which will be explained in the text. First, the right B-model for the equivariant mirror symmetry for Hitchin systems should be given by Landau-Ginzburg models. That is the case for Coulomb branches of $3d$ $\mathcal{N} = 4$ theories, as we briefly mentioned when describing the general setting of Aganagic's proposal. Second, the dichotomy between additive and multiplicative varieties is also crucial. Additive varieties should give the A-model side of the mirror, and multiplicative varieties should appear in the Landau-Ginzburg B-model side. In the last section, building on ideas from hypertoric varieties, this will be made more explicit.\\

Other aspects should be considered in a proposal of (equivariant) mirror symmetry for Hitchin systems. In particular, the role of quantization of Hitchin systems is an interesting subject. On the other hand, as it was mentioned at the beginning, mirror symmetry also appears in three and four dimensions. Therefore, studying the relations between two-dimensional mirror symmetry and three and four-dimensional mirror symmetry seems relevant. Along these lines, the recent work \cite{conan} looks quite important. The relevant discussions regarding homological mirror symmetry, which is part of Aganagic's proposal, are also absent in this text. We plan to discuss some of these topics elsewhere. On the other side, relations between our picture and the so-called Moore-Tachikawa conjecture should also be considered\footnote{We would like to thank David Ben-Zvi for pointing out the relations with the Moore-Tachikawa conjecture. As he explained to us, while Hitchin spaces appear as a Coulomb branch for a class $\mathcal{S}$ theory, its Higgs branch is a more exotic object that is the subject of the Moore-Tachikawa conjecture.}. \\

Finally, we must mention that another important topic, which is not covered in this note, is the relation between our proposal and the geometric Langlands and mirror symmetry for Hitchin systems proposed by Kapustin and Witten (see \cite{kw, wit}). This is really an interesting topic, and its absence could be seen as a big flaw in this text. However, for the initial proposal that we are presenting here, this discussion is not essential, and we can postpone it for other work. However, we are conscious that if we want to really have a complete picture of the mirror symmetry of Hitchin systems, we can not ignore Kapustin and Witten's work. We also plan to pursue this line of thought elsewhere.  

\section{\textbf{Coulomb branches of $3d$ $\mathcal{N}=4$ theories}}

Here we will mainly follow \cite{brav}. Also see \cite{b2,b3,Nak1}. In physics, there is a notion of a $3d$ $\mathcal{N}=4$ supersymmetric quantum field theory. The moduli space $\mathcal{M}$ of vacua of these theories can be a little bit complicated. However, one interesting feature in this case is that $\mathcal{M}$ has two pieces called the Higgs and the Coulomb branches, denoted by $\mathcal{M}_H$ and  $\mathcal{M}_C$ respectively, which are easier to consider. Both $\mathcal{M}_H$ and  $\mathcal{M}_C$ are supposed to be holomorphic symplectic varieties and, in fact, hyperk\"{a}hler. \\

To understand what $\mathcal{M}_H$ and $\mathcal{M}_C$ could be for a $3d$ $\mathcal{N}=4$ theory, it is helpful to consider the following situation: \\

Let $G$ be a complex reductive algebraic group and $\mathrm{M}$ a (symplectic) representation of $G$. That is, $\mathrm{M}$ is a vector space with a symplectic form $\omega$. In addition, assume that the action of $G$ is Hamiltonian. Then, to the pair $(G,\mathrm{M})$ can be associated a theory $\mathfrak{T}(G,\mathrm{M})$. The Higgs branch is expected to be the Hamiltonian reduction of $\mathrm{M}$ with respect to $G$. On the other hand, the Coulomb branch is more complicated to define. The mathematical definition of Coulomb branches of a $3d$ $\mathcal{N}=4$ gauge theory has been undertaken in a series of papers \cite{b2,b3,Nak1} by Braverman, Finkelberg and Nakajima under the assumption that the corresponding theory is of the so-called cotangent type.\\

Now, we will summarize the definition of a Coulomb branch in the spirit of Braverman-Finkelberg-Nakajima. Let us consider a complex connected reductive group $G$ and a finite-dimensional representation $\mathrm{N}$ of $G$. The symplectic representation $\mathrm{M}$ is given by $\mathrm{N} \bigoplus \mathrm{N}^*$. This is the cotangent type condition. Under this assumption, we will consider the pair $(G,\mathrm{N})$ and the construction will be given in terms of $\mathrm{N}$ instead of $\mathrm{M}$. \\

Let $\mathcal{R}_{G, \mathrm{N}}$ the moduli spaces of triples $(\mathcal{P}, \sigma, s)$, where $\mathcal{P}$ is a line bundle on the formal disc $D = \operatorname{Spec}\mathbb{C}[\![z]\!]$, $\sigma$ is trivialization of $\mathcal{P}$ on the formal punctured disc $D^* =\operatorname{Spec}\mathbb{C}(\!(z)\!)$ and $s$ is a section of the associated vector bundle $\mathcal{P}_{\mathrm{triv}} \times_G \mathrm{N}$ on $D^*$ such that $s$ extends to a regular section of $\mathcal{P}_{\mathrm{triv}} \times_G \mathrm{N}$ on $D$ and $\sigma(s)$ extends to a regular section of $\mathcal{P}\times_G \mathrm{N}$ on $D$. \\

Let us denote $\mathbb{C}[\![z]\!]$ by $\mathcal{O}$ and let $G_{\mathcal{O}} = G[\![z]\!]$ be the group of $\mathcal{O}$-valued points of $G$. Then, the $G_{\mathcal{O}}$-equivariant Borel-Moore homology for the moduli space $\mathcal{R}_{G, \mathrm{N}}$, denoted by $H_{\bullet}^{G_{\mathcal{O}}}(\mathcal{R}_{G, \mathrm{N}}$), can be defined (see for instance section 5 of \cite{brav}). This is a commutative, finitely generated and integral algebra. Thus, we have that
the spectrum $\operatorname{Spec}H_{\bullet}^{G_{\mathcal{O}}}(\mathcal{R}_{G, \mathrm{N}})$ of the algebra $H_{\bullet}^{G_{\mathcal{O}}}(\mathcal{R}_{G, \mathrm{N}}$) is defined to be the \emph{Coulomb branch} $\mathcal{M}_C(G,\mathrm{N})$. \\

For our purposes in the last section, it is important to note that if we take $T \subset G$, the maximal torus of $G$, $\mathfrak{t} \subset \mathfrak{g}$ its Lie algebra and $W = \mathrm{N}_G(T)/T$ the corresponding Weyl group, then the  equivariant cohomology $H^{\bullet}_{G_{\mathcal{O}}}(pt)= \mathbb{C}[\mathfrak{t}/W]$ is a subalgebra of  $H_{\bullet}^{G_{\mathcal{O}}}(\mathcal{R}_{G, \mathrm{N}}$) and we have the projection $\pi: \mathcal{M}_C(G,\mathrm{N}) \rightarrow \mathfrak{t}/W$, which is a \emph{integrable system map} whose fibers are Lagrangians. The general fiber is isomorphic to $T^{\vee}$ (dual torus). This is according to the expectation that $\mathcal{M}_C(G, \mathrm{N})$ is approximated by $T^*T^{\vee}/W$. In fact, $\mathcal{M}_C(G, \mathrm{N})$ is birrationally isomorphic to $T^*T^{\vee}/W$. \\

There is an interesting duality between the Higgs branch $\mathcal{M}_H({\mathfrak{T}})$ and Coulomb branch $\mathcal{M}_C({\mathfrak{T}})$ of a $3d$ $\mathcal{N}= 4$ theory $\mathfrak{T}$ and the Higgs branch $\mathcal{M}_H({\mathfrak{T}^*})$ and Coulomb branch $\mathcal{M}_C({\mathfrak{T}^*})$ of another theory $\mathfrak{T}^*$. Physicists expect that for a theory $\mathfrak{T}$ there should exist a theory $\mathfrak{T}^*$ such that $\mathcal{M}_H({\mathfrak{T}^*}) = \mathcal{M}_C({\mathfrak{T}})$ and $\mathcal{M}_H({\mathfrak{T}}) = \mathcal{M}_C({\mathfrak{T}^*})$. This phenomenon is known in the physics literature as \emph{three-dimensional mirror symmetry} (see \cite{phy1}) and in mathematics as \emph{symplectic duality} (see \cite{math1, math2}). 

\subsection{Coulomb branches for quiver gauge theories of type ADE} There are two settings where $\mathcal{M}_C(G, \mathrm{N})$ could be defined for a quiver gauge theory of type ADE, namely the unframed one and the framed one. Here, we will only deal with the framed setting. \\

Let $Q=(Q_0,Q_1)$ be an ADE quiver, where $Q_0$ is the set of vertices and $Q_1$ is the set of arrows. Consider $V = \bigoplus V_i$ a $Q_0$-graded vector space. The corresponding quiver gauge theory is given by the gauge group $\mathrm{GL}(V) = \prod \mathrm{GL}(V_i)$ and the representation $\mathrm{N} = \bigoplus_{\alpha \in Q_1} \mathrm{Hom}(V_{t(\alpha)}, V_{h(\alpha)})$, where $t,h: Q_1 \rightarrow Q_0$ are the tail and the head of an arrow $\alpha$.  Since we have been dealing with the framed case, we also need to consider a $Q_0$-graded vector space $W = \bigoplus W_i$ and add $\bigoplus \mathrm{Hom}(W_i, V_i)$ to $\mathrm{N}$. \\

It was known in the physics literature \cite{witten} (at least in type A) that $\mathcal{M}_C(G, \mathrm{N})$ is a moduli space of singular $G_c$-monopoles on $\mathbb{R}^3$, where $G_c$ is the maximal compact subgroup of $G$. \\

On the other hand, Nakajima conjectured (see \cite{Nak1}) that $\mathcal{M}_C(G, \mathrm{N})$ is a framed moduli space of $S^1$-equivariant instantons on $\mathbb{R}^4$, provided that the coweight $\mu : S^1 \rightarrow G_c$ is dominant, where $mu$ is given by the dimension vectors for ordinary vertices, attached at 0, $\infty$ of $\mathbb{R}^3$. It is important to remark (see \cite{b3}) that there is a subtle difference between $S^1$-equivariant instantons and singular monopoles. As we pointed out above, $S^1$-equivariant instantons only make sense when $\mu$ is dominant. However, they are expected to be isomorphic as complex manifolds.\\

Let $\Lambda$ be the coweight lattice and let $\Lambda^+ \subset \Lambda$ be the submonoid of dominant coweights. Then, for $\lambda \in \Lambda^+$ and $\mu \in \Lambda$, the Coulomb branches for a type ADE quiver can also be identified with a generalized slice $\overline{\mathcal{W}}_{\mu}^{\lambda}$ affine Grassmannian of the corresponding group. For an arbitrary $\mu$, i.e. not necessarily dominant, $\overline{\mathcal{W}}_{\mu}^{\lambda}$ is the moduli space defined by the data (see \cite{FM} for details): \\ 

1. A $G$-bundle $\mathcal{P}$ on $\mathbb{P}^1$. \\
2. A trivialization $\sigma: \mathcal{P}_{\mathrm{triv}} \vert_{\mathbb{P}^1 \setminus \{0 \}} 
\rightarrow \mathcal{P}\vert_{\mathbb{P}^1 \setminus \{0 \}}$ having a pole of degree $\leqslant \lambda$ at $0 \in \mathbb{P}^1$. \\
3. A $\mathrm{B}$-structure $\phi$ on $\mathcal{P}$ of degree $w_0\mu$ with the fiber $\mathbb{B}_{-} \subset G$ at $\infty \in \mathbb{P}^1$, where $\mathbb{B}_{-}$ is the Borel subgroup opposite to $B$ and $w_0 \in W$ is the longest element. \\

The slice $\overline{\mathcal{W}}_{\mu}^{\lambda}$ is nonempty if and only if $\mu \leq \lambda$.\\

If the quiver $Q$ is an affine quiver of type $\tilde{A}$ $\tilde{D}$ $\tilde{E}$, it is expected that $\mathcal{M}_C$ is isomorphic to the Uhlenbeck partial compactification of the moduli space of $G$-bundles on $\mathbb{P}^2$ trivialized at $\mathbb{P}^1_{\infty}$, of second Chern class $d$. When $G = SL_n(\mathbb{C})$ this was proved in \cite{NT}.

\section{\textbf{Aganagic's proposal of equivariant mirror symmetry}}\label{agan}

Motivated by the problem of categorifying quantum knot invariants, Aganagic \cite{ag1,ag2} (see also \cite{ag}) proposed an approach to mirror symmetry for certain Coulomb branches of a $3d$ $\mathcal{N}=4$ quiver gauge theory, taking into account the action of a torus on such Coulomb branches. \\ 

Let $\mathcal{X}$ be a Coulomb branch of an $3d$ $\mathcal{N}=4$ quiver gauge theory with quiver $Q$, the Dynkin diagram of $\mathfrak{g}$, where $\mathfrak{g}$ is the Lie algebra of a group $G$. According to the discussion of the previous section $\mathcal{X}$ also has an interpretation as the moduli space of singular $G$-monopoles on $\mathbb{R}^3 = \mathbb{R} \times \mathbb{C}$ or as a resolution of a slice in the affine Grassmannian. For us, the interpretation as Coulomb branches or as moduli spaces of certain monopoles will be the most relevant due to the relation with physics. \\

As it was mentioned in the introduction, Aganagic's proposal for the equivariant mirror of $\mathcal{X}$ works at two levels. One level is the \emph{upstairs mirror} of $\mathcal{X}$, and the other one is the \emph{downstairs mirror} of its core $X$. The equivariant mirror of $\mathcal{X}$ would be the \emph{downstairs} Landau-Ginzburg model mirror to $X$, whose superpotential $W$ mirrors the action of a torus $T$ on $\mathcal{X}$. Let us see this in more detail. 

\subsection{Upstairs mirror} Let us construct the two-dimensional mirror for $\mathcal{X}$ as a holomorphic symplectic manifold. To do this, it is convenient to consider the Coulomb branch $\mathcal{Y}^\vee$ of a four-dimensional gauge theory based on the same quiver $Q$ and compactified on a circle of radius $R^\vee$. \\

The first thing to note is that $\mathcal{Y}^\vee$ is also holomorphic symplectic. The circle compactification reduces the dimension from four to three. By further compactifying on a circle of radius $R$, we obtain a two-dimensional sigma model with target space $\mathcal{Y}^\vee$. Now, since $\mathcal{Y}^\vee$ is also hyperk\"aler, its mirror $\mathcal{Y}$ can be obtained by hyperk\"ahler rotation.\\

An important remark, due to Aganagic, is that $\mathcal{Y}$ can be interpreted as the multiplicative version of $\mathcal{X}$, and it is part of the mirror of $\mathcal{X}$. It is not exactly the mirror of $\mathcal{X}$, because the precise comparison between $\mathcal{X}$ and $\mathcal{Y}^\vee$ leads to the consideration of a superpotential $\mathcal{W}$. \\

In the limit when the radius $R^\vee$ goes to zero, $\mathcal{Y}^\vee$ becomes the Coulomb branch $\mathcal{X}$ of the three-dimensional gauge theory with quiver $Q$. More precisely, we have that, as holomorphic symplectic manifolds,

\begin{equation*}
\mathcal{Y}^\vee = \mathcal{X} \setminus \mathrm{D}^\vee
\end{equation*} 

where $\mathrm{D}^\vee$ is the union of some divisors $\mathrm{D}_a^\vee$. Each $\mathrm{D}_a^\vee$ is the locus in $\mathcal{X}$ where $\mathrm{det}(\Phi_a) = 0$ and $\Phi_a$ is the complex scalar valued in the adjoint representation of the gauge group. Note that both $\mathcal{X}$ and $\mathcal{Y}^\vee$ are both Coulomb branches but for a three-dimensional theory and four-dimensional theory, respectively. \\

We can summarize the above discussion in the following diagram 

\[
\begin{tikzcd}[column sep=4em, row sep=4em]
\mathcal{Y}^\vee = \mathcal{X} \setminus \mathrm{D}^\vee
  \arrow[r, rightarrow, "R^\vee \rightarrow 0"] 
  \arrow[d, "mirror"] 
& \mathcal{X} \\
\mathcal{Y}
\end{tikzcd}
\]

Now, to obtain the actual mirror of $\mathcal{X}$, it is needed to understand the mirror to a partially compactifying $\mathcal{Y}^\vee$ by adding in the divisors $\mathrm{D}_a^\vee$. The addition of divisors corresponds, by mirror symmetry, to turning on a superpotential $\mathcal{W}$ on $\mathcal{Y}$, so the two-dimensional mirror of $\mathcal{X}$ is a Landau-Ginzburg model $(\mathcal{Y},\mathcal{W})$. This is the \emph {upstairs mirror}. \\

We should note that, so far, we have not considered the action of a torus in the picture. Before discussing this, let us note that this kind of Landau-Ginzburg models $(\mathcal{Y},\mathcal{W})$ also appear in the recent work of Khan and Moore \cite{moore} and in the work of Kontsevich and Soibelman \cite{exponential}. They also appear (implicitly) in the work of the first author with Chen and Romo about an equivariant version of the Dubrovin conjecture \cite{cruz}.

\subsection{Downstairs mirror}

We have that $\mathcal{X}$ has a $\mathbb{C}^{\times}_q$-symmetry which scales its holomorphic symplectic form. The invariant locus $X$ sits inside $\mathcal{X}$ as a holomorphic Lagrangian and it is called the \emph{core} of $\mathcal{X}$. In addition, $\mathcal{X}$ has a larger torus of symmetries, namely $T= (\mathbb{C}^{\times})^{\mathrm{rk}\mathfrak{g}} \times \mathbb{C}^{\times}_q$, where $\mathfrak{g}$ is the Lie algebra of $G$. We want to consider the action of $T$ for constructing the `\emph{downstairs mirror}. \\

The \emph{downstairs mirror} is the Landau-Ginzburg model $(Y,W)$ which is the mirror of the core $X$. Roughly speaking, it is related to the the \emph{upstairs mirror} in the following sense: the superpotential $W$ is a multivalued holomorphic function which mirror the equivariant $T$-action on $\mathcal{X}$ and $\mathcal{Y}$ fibers over $Y$ with holomorphic Lagrangians $(\mathbb{C}^{\times})^m$-fibers, where $m$ is the dimension of $X$, i.e., half of the dimension of $\mathcal{X}$. \\

Thus, thanks to the fibration above, there is a relation between the superpotential $\mathcal{W}$ and $W$. Consider a local parametrization of $\mathcal{Y}$ in terms of coordinates $y_{a,\alpha}$, and local coordinates $z_{a,\alpha}$ valid for large $y$. \\

The form of the superpotential $\mathcal{W}$ on $\mathcal{Y}$ must satisfy: \\

1. Restricting to $\partial_z \mathcal{W} = 0$, we recover the superpotential $W$. \\

2. The terms of $\mathcal{W}$ can be written in terms of globally defined coordinates on $\mathcal{Y}$. \\

Using these constrains, Aganagic conjectures the form of $\mathcal{W}$ as (see \cite{ag2}))

\begin{equation*}
\mathcal{W} = \sum_a \sum_{\alpha} (\lambda_a \operatorname{log}(y_{a,\alpha}) - \lambda_0 \operatorname{log}(z_{a,\alpha}) + z_{a,\alpha}^{-1}f_{a,\alpha}(y))
\end{equation*}
where, $f_{a,\alpha}(y)$ are certain holomorphic functions. See \cite{ag2} for details about these functions.\\ 

This $\mathcal{W}$ has some additional terms due to the action of $T$, so it is slightly different than the $\mathcal{W}$ defined in the \emph{upstairs mirror} above. However, abusing the notation, we write both as the same. \\

Now, we can improve the mirror diagram in the introduction to obtain 

\[
\begin{tikzcd}[column sep=6em, row sep=4em]
(\mathcal{X},T)
  \arrow[r, rightarrow, "upstairs-mirror"] 
  \arrow[dr, "equiv.mirror" description] 
& (\mathcal{Y},\mathcal{W}) \arrow[d, rightarrow, "\partial_z \mathcal{W} = 0 \rightarrow W"] \\
X \arrow[r, rightarrow, "downstairs-mirror"']
 \arrow[u, hook]  
& (Y,W)
\end{tikzcd}
\] \\

Thus, we obtain the following Theorem

\begin{teorema}\label{teor1}[Aganagic]
The Landau-Ginzburg model $(Y,W)$ is the equivariant mirror of $\mathcal{X}$, in the sense that the A-model of $\mathcal{X}$, computed by Gromov-Witten theory of $\mathcal{X}$ working equivariantly with respect to the action of $T$ coincides with the B-model of $(Y,W)$, even though $\mathcal{X}$ and $Y$ do not have the same dimension. In fact, the correspondence holds in all genera (not only in genus zero).
\end{teorema}

The last sentence comes from the fact that the genus zero theory is semisimple, and applying Givental-Teleman reconstruction, one can extend the correspondence to all genera.

\section{\textbf{Hitchin systems and class $\mathcal{S}$ theories}}\label{classS}

The question is whether the kind of equivariant mirror symmetry described in the last section can be extended beyond the situation studied by Aganagic. The motivation of this note is to explore one possible scenario for such an extension. For this, we need to consider Hitchin systems, from the mathematical side, and class $\mathcal{S}$ theories from the physical side.  

\subsection{Hitchin systems} The literature about Hitchin systems is huge. Here, we only recall some basic definitions. For a more detailed and recent presentation, see \cite{Etingof}. \\

In \cite{hitchin87}, Hitchin studied the moduli space of solutions to certain self-duality equations on a Riemann surface $C$. Let us denote it by $\mathcal{M}_{\mathrm{Hit}}$. For a complex reductive group $G$, $\mathcal{M}_{\mathrm{Hit}}$ carries a (canonical) hyperk\"ahler metric with complex structures $I, J, K$. \\

Let $E$ be a $G$-bundle and $\phi \in \mathrm{H}^0(C, \mathrm{ad}(E) \otimes K)$, where $K$ is the canonical bundle of $C$ and $\mathrm{ad}(E)$ is the adjoint bundle. This $\phi$ is called Higgs field. Consider the moduli space of pairs $(E,\phi)$ (moduli space of Higgs bundles), and denote it by $\mathcal{M}_{\mathrm{Dol}}$. Taking the complex structure $I$ in $\mathcal{M}_{\mathrm{Hit}}$, we have that $(\mathcal{M}_{\mathrm{Hit}}, I) \simeq \mathcal{M}_{\mathrm{Dol}}$. For the other complex structures $J$ and $K$ we have that $(\mathcal{M}_{\mathrm{Hit}}, J) \simeq (\mathcal{M}_{\mathrm{Hit}}, K) \simeq \mathcal{M}_{\mathrm{DR}}$, where $\mathcal{M}_{\mathrm{DR}}$ is the moduli space of flat $G$-connections on $C$. \\

One the most important features related to these moduli spaces is the existence of the so-called Hitchin map (or Hitchin fibration)
\begin{equation*}
h: \mathcal{M}_{\mathrm{Dol}} \rightarrow \mathcal{A} = \bigoplus_{i=1}^n \mathrm{H}^0(C, K^i)
\end{equation*}
given by $(E,\phi) \mapsto \mathrm{det}(\phi - \lambda I)$. The dimension of $\mathcal{A}$ is half the dimension of $\mathcal{M}_{\mathrm{Dol}}$, and the fibers are Lagrangian subvarieties. Hitchin proved that $h$ is an algebraically completely integrable system. For this reason, it is known as \emph{Hitchin system}. We want to remark that sometimes, we may require some punctures on $C$. This will be the case for the Hitchin systems coming from class $\mathcal{S}$ theories, which we will discuss below. \\ 

There are two recent developments in the theory of Hitchin systems that are relevant for our purposes. The first one is the work of Hausel and Hitchin \cite{haus2} (see also \cite{haus1}) on an enhanced mirror symmetry for Hitchin systems and the second one is the work of Elliot and Pestun \cite{elliot} (see also \cite{elliot2, nekrasov-pestun}) on the relation between a multiplicative version of Hitchin systems and supersymmetric gauge theories. 

\subsubsection{\textbf{Enhanced mirror symmetry}} 

Based on the previous work \cite{Hauselmirror}, Hausel and Hitchin propose a setting for considering mirror symmetry, at the categorical level, for Hitchin systems. They were interested in considering the $\mathbb{C}^{\times}$-action on $\mathcal{M}_{\mathrm{Hit}}$ given by $(E, \phi) \mapsto (E, \lambda \phi )$, the scalar multiplication of the Higgs field. Note that the Hitchin map is also $\mathbb{C}^{\times}$-equivariant if $\mathbb{C}^{\times}$ acts on $\mathrm{H}^0(C, K^i)$ with weight $i$. \\

With respect to the $\mathbb{C}^{\times}$ action defined above, the core $\mathfrak{C}$ of $\mathcal{M}_{\mathrm{Hit}}$ can be defined. It was pointed out in \cite{haus2} that the core can be seen as the reduced scheme of the nilpotent cone $h^{-1}(0)$. It was also observed that the core is a Lagrangian subvariety. The action and the defined core lead Hausel and Hitchin to consider an enhanced mirror symmetry for Hitchin systems.\\ 

As far as we understand, the enhancement is given by the consideration of the torus action in the mirror picture. There is a tempting analogy with Aganagic's proposal concerning the \emph{downstairs mirror}. We will return to this point later. Another interesting thing is the study of a $\mathbb{C}^{\times}$-equivariant extension of the Donagi-Pantev correspondence $D^b(\mathcal{M}_{\mathrm{Dol}}) \simeq D^b(\mathcal{M^\vee}_{\mathrm{Dol}})$, which is a classical limit of the geometric Langlands correspondence. We will not discuss this here, but the interested reader could see \cite{haus2}.\\

Formally, we have the following comparison between Aganagic's picture and Hausel-Hitchin's picture

\begin{table}[h]
\centering
{%
\begin{tabular}{ |c|c| }
 \hline
 \textbf{Coulomb branch of a $3d$ $\mathcal{N}=4$ theory} & \textbf{Hitchin system}  \\ \hline
 $\mathbb{C}^{\times}_{\mathfrak{q}}$ symmetry scaling the holomorphic form  & $\mathbb{C}^{\times}$ action \\  
 core $X$ &  core $\mathfrak{C}$ defined as \\ 
       &           the reduced nilpotent cone             \\
 LG model $(Y,W)$ mirror to X& mirror to  $\mathfrak{C}$ ?  \\
 \hline
\end{tabular}%
}

\end{table} 

\subsubsection{\textbf{Multiplicative Hitchin systems}}  We will follow the developments in \cite{elliot}, but the moduli spaces of multiplicative Higgs bundles appear in the literature earlier under the name of moduli of $G$-pairs, see for instance \cite{hurtubise, hurtubise1} for the case when $G$ is a general reductive group. \\

A multiplicative Higgs bundle on $C$ (Riemann surface) is a $G$-bundle $P$ on $C$ with a meromorphic automorphism $g$, this means than in the usual definition we have to replace the adjoint bundle $\mathrm{ad}(P)$ by the group valued adjoint bundle $\mathrm{Ad}(P)$. The moduli space of these objects admits a version of the Hitchin map. It is important to remark that even though the multiplicative map makes sense for any curve, it only defines a completely integrable system in the case when the curve $C$ is Calabi-Yau, i.e. $C$ is one of the following curves: $\mathbb{C}$, $\mathbb{C}^{\times}$ or an elliptic curve. This is particularly important, since in the usual definition of the Higgs fields, we take a section of the adjoint bundle twisted by the canonical bundle. It is not clear at all what the replacement for this twist is in the multiplicative case; however, in the Calabi-Yau case, such a twist is trivial, so we do not need to worry about it. \\

For the Calabi-Yau case, the resulting moduli space (under some conditions, see \cite{elliot} for details) is hyperk\"ahler. We only note that if $C = \mathbb{C}$, $\mathbb{C}^{\times}$, it is needed to take a specific boundary conditions at $\infty$.\\

From a physics point of view, there are two main sources for multiplicative Hitchin systems. In \cite{elliot}, they appear as the moduli spaces of solutions to equations of motion in certain 5-dimensional supersymmetric gauge theories. On the other hand, in \cite{nekrasov-pestun} they appear as the Seiberg-Witten integrable system of a $4d$ $\mathcal{N}=2$ theory. The former approach is motivated by the idea of finding a multiplicative version of the geometric Langlands conjecture. Although this could be relevant for us, the approach in \cite{nekrasov-pestun} is better for our proposal, since they identify integrable systems underlying the special geometry of the $\mathcal{N}=2$ theories with Higgs configurations where the gauge group $G$ corresponds to a quiver diagram. On the other hand, in \emph{loc.cit.}, they point out that having a complete description of the Seiberg-Witten curves and algebraic integrable systems for the $\mathcal{N}= 2$ ADE quiver theories would be interesting to further investigate this ADE quiver class along the lines of class $\mathcal{S}$ theories.

\subsection{Class $\mathcal{S}$ theories} \label{classS} We will follow the motto expressed by Neitzke in \cite{Neitzke2014}, namely, \\

``\emph{the relation between $\mathcal{N}=2$ theories and Hitchin systems is that the Hitchin system arises as the moduli space of the $\mathcal{N}=2$ theory compactified on a circle}". \\

It is also pointed out in \emph{loc.cit.} that the $\mathcal{N}=2$ theories mostly related to Hitchin systems are the so-called \emph{class $\mathcal{S}$ theories}. So, What is a class $\mathcal{S}$ theory?  Let us start answering this question by considering a general picture. \\

Consider a (2,0) superconformal field theory in 6 dimensions of type ADE, and compactifying it on an appropriate (punctured) Riemann surface $\Sigma$. This compactification procedure produces a four-dimensional $\mathcal{N}=2$ field theory. See \cite{GMN} for details. We are interested in a further compactification on a circle $\mathrm{S}^1$ of radius $R$ of this four-dimensional theory. After these two compactifications, at low energies, we end with a three-dimensional sigma model with a hyperk\"ahler target space $\mathcal{M}$.\\

We would like to remark that there is also a possibility to reverse the compactification and obtain the same three-dimensional sigma model. Starting with the six-dimensional theory and compactifying it on a circle $\mathrm{S}^1$ of radius $R$, we obtain a five-dimensional supersymmetric Yang-Mills theory. Then, compactifying this theory on $\Sigma$, with an appropriate topological twist, we get the three-dimensional theory. See \cite{Neitzke2014} for more details.\\

The class $\mathcal{S}$ theories are the three-dimensional theories obtained from one of the compactification procedures described above. This is equivalent to saying that the three-dimensional theory is obtained from a six-dimensional one from a compatification on $\Sigma \times \mathrm{S}^1$. \\

Now, we can ask how Hitchin systems appear in this setting. Let us consider the $4d$ $\mathcal{N} = 2$ supersymmetric theory, which we denote by $T$, and let $\mathfrak{B}$ be its Coulomb branch. $\mathfrak{B}$ consists of an open regular locus $\mathfrak{B}_{\mathrm{reg}}$ and a discriminant locus $\mathfrak{B}_{\mathrm{sing}}$. It is known (see \cite{Neitzke2014}) that there is no single Lagrangian describing the theory globally, but there is a single geometric object from which all local Lagrangians can be described. That object is a complex integrable system associated with the theory $T$. More concretely, there is a holomorphic symplectic manifold $\mathcal{I}'$ with a projection $\pi: \mathcal{I}' \rightarrow \mathfrak{B}_{\mathrm{reg}}$ such that $\mathcal{I}'_u = \pi^{-1}(u)$ are compact complex Lagrangian tori. We want to remark that this can be extended to the whole $\mathfrak{B}$, but for our discussion it is enough to take $\mathfrak{B}_{\mathrm{reg}}$. \\

The three-dimensional theory obtained after compactifying on $S^1$ is a sigma model whose target $\mathcal{M}$ is related to the complex integrable system $\mathcal{I}'$. The constrains of $\mathcal{N}=4$ supersymmetry in dimension three force that the metric on $\mathcal{M}$ (an actual sigma model needs a Riemannian metric) must to be hyperk\"ahler, i.e., that it carries a family of complex structures $J_{\xi}$, $\xi \in \mathbb{CP}^1$, as well as corresponding holomorphic symplectic forms $\omega_{\xi}$. One of these complex structures is distinguished, let us denote it $J_0$, and when considering as a holomorphic symplectic manifold in $J_0$, $\mathcal{M}$ is identical to $\mathcal{I}'$. For class $\mathcal{S}$ theories, the complex integrable systems we obtain are Hitchin systems. \\

The following diagram summarized the construction
\[
\begin{tikzcd}[column sep=6em, row sep=4em]
\mathcal{M}
& \mathcal{I}' 
\arrow[l, rightarrow, "3d-limit"]
\arrow[d, rightarrow, "Hit.fibration"] \\
& \mathcal{B}_{\mathrm{reg}}
\end{tikzcd}
\] 

Formally, we have the following comparison between Aganagic's picture and the Hitchin system associated to class $S$ theories

\begin{table}[h]
\centering
{%
\begin{tabular}{ |c|c| }
 \hline
 \textbf{Coulomb branch of a $3d$ $\mathcal{N}=4$ theory} & \textbf{Hitchin system asso. to class $\mathcal{S}$}  \\ \hline
 $\mathcal{Y}^\vee$ Coulomb branch of a $4d$ theory  & $\mathfrak{B}_{\mathrm{reg}}$ Coulomb branch of a $4d$ theory \\  
 $\mathcal{X}$ a 3-dim limit &  $\mathcal{M}$ a 3-dim limit \\ 
 $\mathcal{Y}^\vee$ mirror to $\mathcal{Y}$ & $\mathcal{I}'$ self-mirror? \\
 $\mathcal{Y} \rightarrow Y$ fibration & $\mathcal{I}' \rightarrow \mathfrak{B}_{\mathrm{reg}}$ fibration  \\
 $\mathcal{Y}^\vee = \mathcal{X} \setminus D^\vee$ as complex manifolds &  $\mathcal{I}' = \mathcal{M} \setminus D$ for a \\
    & complex structure $J_{\xi}$, with $\xi \neq 0$ and a divisor $D$? \\
 
 \hline
\end{tabular}%
}

\end{table} 

Finally, we want to mention the work of Shan, Xie and Yan \cite{shan}. Motivated by three-dimensional mirror symmetry for $3d$ $\mathcal{N}=4$ superconformal field theories they propose a mirror symmetry picture for a $4d$ $N=2$ superconformal field theories compactified on $\mathrm{S}^1$ with finite radius $R$, in terms of a duality (they call it algebra/geometry duality) between certain vertex operator algebras associated to the $4d$ theory and the Coulomb branch. They focus on the case when the $4d$ theory is of class $\mathcal{S}$ and provide a new perspective for mirror symmetry of Hitchin systems in terms of their algebra/geometry duality.

\section{\textbf{First steps to an equivariant mirror symmetry for Hitchin systems}}\label{proposal}

This section would be highly speculative. It is based on some formal similarities between the contexts described in the previous sections, and not on concrete evidence. To make our proposals more plausible, we want to start discussing the work of Gammage, McBreen, and Webster in the context of hypertoric varieties, which, at least philosophically, is in the same spirit as this note. \\

In \cite{GMW}, mirror symmetry for multiplicative hypertoric varieties was studied. The authors point out that additive hypertoric varieties appear as Coulomb branches of $3d$ $\mathcal{N}=4$ gauge theories, like the variety $\mathcal{X}$ in Aganagic's work. On the other hand, multiplicative hypertoric varieties appear as Seiberg-Witten systems governing a $4d$ $\mathcal{N}=2$ theories, like $\mathcal{Y^\vee} $ and $\mathcal{Y}$ in Aganagic's work and the Hitchin system associated to a class $\mathcal{S}$ theories. \\

In particular, Gammage, McBreen and Webster conjecture

\begin{conjetura}[\cite{GMW}] \label{con} Let $\mathfrak{U}$ be a multiplicative hypertoric variety. Then there exists an additive variety $\mathfrak{M}$ a divisor $\mathfrak{D} \in \mathfrak{M}$ and an isomorphism $\mathfrak{U} \simeq \mathfrak{M} \setminus \mathfrak{D}$. 
\end{conjetura}

We should remark that the authors of \cite{GMW} do not expect the conjecture to hold beyond the hypertoric case. On the contrary, as we will establish below, we believe that there is room for considering an \emph{extension of multiplicative and additive varieties} in which a similar conjecture should hold. \\

Related to conjecture \ref{con}, Lau and Zheng \cite{lau} proved the following theorem

\begin{teorema} \label{teor} Let $\mathfrak{M}$ be an additive hypertoric varieties, then $\mathfrak{M}$ is SYZ mirror to a Landau-Ginzburg model $(\mathfrak{U}^\vee, W: \mathfrak{U}^\vee \rightarrow \mathbb{C})$ whose underlying space $\mathfrak{U}^\vee$ is a multiplicative hypertoric variety. In particular, the complement in $\mathfrak{M}$ of a divisor $\mathfrak{D}$ is SYZ mirror to the multiplicative hypertoric variety $\mathfrak{U}^\vee$
\end{teorema}

From the discussion in section \ref{agan}, it is clear that taking $\mathfrak{M}$ as the Coulomb branch $\mathcal{X}$ and $\mathfrak{U}^\vee$ as $\mathcal{Y}^\vee$ (or its mirror $\mathcal{Y}$), conjecture \ref{con} holds. Also, the \emph{upstairs mirror} (without the consideration of the torus action) is a consequence of Theorem \ref{teor}. \\

It is natural to think that if we incorporate a torus action on  $\mathfrak{M}$, a core $M$ could be defined. It is also plausible to expect a fibration  $\mathfrak{U}^\vee \rightarrow U^\vee$, such that $M$ and $(U^\vee, W)$ are mirror partners. Thus, we would like to propose:

\begin{proposal} For hypertoric varieties in conjecture \ref{con}, the \emph{upstairs mirror} and the \emph{downstairs mirror}, with respect to the action of a torus $T$, in the sense of the discussion of section \ref{agan}, could be defined such that a version of Theorem \ref{teor1} holds.
\end{proposal}

However, we would like to go beyond the case of hypertoric varieties. The discussion in the sub-section \ref{classS} motivates us to consider the holomorphic symplectic variety $\mathcal{I}'$ as a \emph{multiplicative variety} and the three-dimensional limit $\mathcal{M}$ as an additive variety. We saw that in the particular choice $J_0$ of a complex structure, $\mathcal{I}'$ coincides with $\mathcal{M}$. However, we want to propose:

\begin{proposal} There exist a complex structure $J_\xi$ such that, as complex manifolds $\mathcal{M} \simeq \mathcal{I}' \setminus D$, for some divisor $D$ of $\mathcal{M}$. 
\end{proposal}

This is a variant of conjecture \ref{con} for Hitchin systems associated to class $\mathcal{S}$ theories discussed in \ref{classS}. \\

Thus, we want to propose an \emph{upstairs mirror} in the Hitchin systems situation. 

\begin{proposal} The mirror of (the additive variety $\mathcal{M}$ is a Landau-Ginzburg model $\mathcal{W}: \mathcal{I}' \rightarrow \mathbb{C}$. If we consider a torus action on $\mathcal{M}$, some extra terms should be added to $\mathcal{W}$ reflecting the equivariant parameters of the action. 
\end{proposal}

One of the features of the B-model side in Aganagic's proposal is that the LG model comes with a fibration. For the case of the Hitchin systems, we have the Hitchin fibration $\mathcal{I}' \rightarrow \mathfrak{B}_{\mathrm{reg}}$, so we could think of adding the superpotential $\mathcal{W}$ to this fibration defined in such a way that its restriction according to similar rules as in Aganagic's setting yields a superpotential $W$ for $\mathfrak{B}_{\mathrm{reg}}$. In addition, following Hausel-Hitchin, we could define a core $\mathfrak{C}$ in $\mathcal{M}$ (or a variant of it, depending on the torus action we consider). Therefore, we would like to propose a \emph{downstairs mirror} in the Hitchin situation.

\begin{proposal} There exists a torus $T$ and an action of it on $\mathcal{M}$ such that the core $\mathfrak{C}$ with respect to the action of $T$ could be defined. The core can be defined in such a way that it is a mirror to the Landau-Ginzburg model $(\mathfrak{B}_{\mathrm{reg}}, W)$, where the superpotential $W$ mirrors the equivariant action of $T$.
 \end{proposal}
 
One of the things where Aganagic's picture and the Hitchin systems case differ is that the Coulomb branch for the four-dimensional theory is the total space of the fibration in Aganagic's framework, and the base of the fibration in the Hitchin systems situation. This makes it more difficult to construct the \emph{downstairs mirror} in the Hitchin setting. We think that making a more precise proposal for the \emph{downstairs mirror} is an interesting question. What we believe is that a variant for this mirror should be found in such a way that the following diagram holds for the Hitchin systems associated to class $\mathcal{S}$ theories.

\[
\begin{tikzcd}[column sep=6em, row sep=4em]
(\mathcal{M},T)
  \arrow[r, rightarrow, "upstairs-mirror"] 
  \arrow[dr, "equiv.mirror" description] 
& (\mathcal{I}',\mathcal{W}) \arrow[d, rightarrow, "\partial_z \mathcal{W} = 0 \rightarrow W"] \\
 \mathfrak{C} \arrow[r, rightarrow, "downstairs-mirror"']
 \arrow[u, hook]  
& (\mathfrak{B}_{\mathrm{reg}},W)
\end{tikzcd}
\] \\

One interesting case to consider is when the multiplicative variety is given by the multiplicative Hitchin systems in the sense of Elliot and Pestun. As we noted before, in this case, we obtain an actual integrable system only when the genus of the curve $C$ is 0 or 1. Thus, we expect that our proposals hold in this case. In the general case, when the genus $g$ is bigger than 1, we also have a fibration, so a candidate for a B-model could be considered. However, we do not have any idea of how the A-side could be in this situation. This could also be related in an interesting way to geometric Langlands, see \cite{elliot}. \\

There is also another missing piece, and it is the equivariant Gromov-Witten theory for the Hitchin systems. We need this for formulating the equivariant mirror symmetry properly. Due to the formal similarities between the Hitchin systems associated to theories of class $\mathcal{S}$ and Nakajima quiver varieties, we expect that some of the constructions in the literature for the (equivariant) Gromov-Witten theory of Nakajima quiver varieties (see \cite{MO}) could be mimicked in the case of the Hitchin systems. Therefore, we want to end with the following proposal:

\begin{proposal} The (equivariant) genus 0 Gromov-Witten theory for a Hitchin system of class $\mathcal{S}$ (the holomorphic symplectic manifold $\mathcal{I}'$) can be constructed along the same lines than the one for Nakajima quiver varieties. In this case, we can produce a semisimple Frobenius manifold and applying the reconstruction theorem of Givental-Teleman, we can obtain the full genera information. 
\end{proposal}

\section*{Acknowledgements}

I would like to thank Olivia Dumitrescu and Motohico Mulase for interesting and illuminating discussions related to this note. This note started as a joint collaboration with them and I have been benefited from their insights. However, any inaccuracy is my own responsibility. I would also like to thank David Ben-Zvi and Sergey Galkin for having read a draft version of this note and their useful comments. I also would like to thank the anonymous referee for his/her comments and suggestions that substantially improve the presentation and the content of this paper.\\

The initial part of this paper was written during the time I was a CAS-PIFI fellow as visiting professor of the University of Science and Technology of China. I thank the Chinese Academy of Science and the University of Science and Technology of China for the financial support and excellent working conditions. The last part of the paper was written during my visit to the Hausdorff Research Institute for Mathematics (HIM) and the Max Planck Institute for Mathematics (MPIM). I thank both institutions for their support and excellent working conditions. During my stay at HIM and MPIM I was funded by the Deutsche Forschungsgemeinschaft (DFG, German Research Foundation) under Germany's Excellence Strategy – EXC-2047/1 – 390685813.


\end{document}